\documentclass[a4paper,12pt]{article}
\usepackage{amssymb}
\usepackage{amsmath}
\usepackage{graphicx}

\newcommand{\der}{\mathrm{d}}

\newcommand{\F}{\mathcal{F}}

\newcommand{\M}{\mathcal{M}}
\newcommand{\Mn}{\mathbb{E}}
\newcommand{\N}{\mathbb{N}}
\newcommand{\nep}{\mathrm{e}}
\newcommand{\Pb}{\mathbb{P}}

\newcommand{\Q}{\mathbb{Q}}
\newcommand{\Rl}{\mathbb{R}}
\newcommand{\T}{\mathrm{T}}
\newcommand{\Var}{\mathrm{Var}}
\newcommand{\Z}{\mathbb{Z}}
\newtheorem{thm}{Theorem}
\newtheorem{rmk}[thm]{Remark}
\newtheorem{lem}[thm]{Lemma}
\newtheorem{xpl}[thm]{Example}
\newtheorem{defi}[thm]{Definition}
\newtheorem{cor}[thm]{Corollary}
\title{Coarse Ricci curvature for continuous-time Markov processes}
\author{Laurent Veysseire}
\date{}
\begin{document}
\maketitle
\abstract{In this paper, we generalize Ollivier's notion of coarse Ricci curvature for Markov chains to continuous time Markov processes.
We prove Wasserstein contraction and a Lichnerowicz-like spectral gap bound for reversible Markov processes with positive coarse Ricci curvature.}

\section*{Introduction}

In \cite{olli}, Ollivier defines the coarse Ricci curvature for Markov chains on metric spaces, in a discrete time framework.
Here we extend this notion to continuous time Markov processes.
We define the curvatures $\kappa$ and $\bar{\kappa}$ (see Definition \ref{curv}), and prove (Theorem \ref{inf2}) that a control of $\bar{\kappa}$ (or $\kappa$) implies that the Markov process contracts the $W_1$ Wasserstein distance between measures exponentially fast (or that the $W_1$ distance des not explode faster than exponentially, if we have negative curvature).
Note that the definition of the coarse Ricci curvature is local.
It is natural to think that positive curvature gives global contraction, but it is not a trivial consequence.

We also show that the coarse Ricci curvature allows to generalize the Lichnerowicz Theorem (see \cite{lich}) that we recall below.

\begin{thm}[Lichnerowicz]\label{Lic}
Let $(\M,g)$ be a $n$-dimensional Riemannian manifold.
If there exists $K>0$ such that for each $x\in\M$, for each $u\in\T_x\M$, we have $\mathrm{Ric}_x(u,u)\geq K g_x(u,u)$, then the spectral gap $\lambda_1$ of the Laplace operator $\Delta$ acting on $L^2$ satisfies
$$\lambda_1\geq\frac{n}{n-1}K.$$
Here we denote by $\mathrm{Ric}$ the Ricci curvature tensor of $\M$.
\end{thm}

Using the contraction of the Wasserstein distance given by Theorem \ref{inf2}, we can prove we have a spectral gap.
What we can get using coarse Ricci curvature is the following:
\begin{thm}\label{inf} Let $(P^t)$ be the semi-group of a reversible and ergodic Markov process on a Polish space $(E,d)$, admitting a left continuous modification.
Assume that for all $(x,y)$ in $E^2$ with $x\neq y$, the coarse Ricci curvature $\kappa(x,y)$ between $x$ and $y$ (see Definition \ref{curv}) is bounded below by a constant $K>0$.
Also assume that for some (then any) $x_0$, $\int d^2(x,x_0)\der\pi(x)<+\infty$, where $\pi$ is the reversible measure.

Then the operator norm of $P^t$ acting on the space $L^2_0(E,\pi)$ of $0$-mean $L^2(\pi)$ functions is at most $\nep^{-Kt}$.
In the case when the Markov process admits a generator $L$, this means $L$ has a spectral gap $\lambda_1(L)\geq K$.
\end{thm}

This Theorem looks like Theorem 1.9 of \cite{cw94}.
The difference is that the contraction hypothesis was global, and some assumption about the first eigenfunction was required.
Here in Theorem \ref{inf2}, the coarse Ricci curvature is local.

In the special case of diffusion processes on Riemannian manifolds, we can get better lower bounds for the spectral gap, depending on the harmonic mean of the Ricci curvature instead of its infimum, as shown in \cite{vey10} and \cite{veys}.

\section{Coarse Ricci curvature: definition and examples}
The coarse Ricci curvature of a Markov process on a Polish (metric, complete, separable) space $(E,d)$ is defined thanks to the Wasserstein metric, which is based on optimal coupling:
\begin{defi}
The Wasserstein distance between two probability measures is the (possibly infinite) quantity defined by:
$$W_1(\mu,\nu)=\inf_{\xi\in\Phi(\mu,\nu)}\int d(x,y)\der\xi(x,y).$$
Here $\Phi(\mu,\nu)$ is the set of all couplings between $\mu$ and $\nu$, that is, the set of probability measures on $E^2$ whose marginal laws are $\mu$ and $\nu$.
\end{defi}
The duality theorem of Kantorovitch (see \cite{vill}) gives another interpretation of this distance and allows to extend it to finite measures provided they have the same total mass, and makes the triangular inequality for $W_1$ easier to check.
\begin{thm}[Kantorovitch--Rubinstein]We have the following equality:
$$W_1(\mu,\nu)=\sup_{f\,\mathrm{bounded, }1\mathrm{-Lipschitz}}\int f\der(\mu-\nu).$$
\end{thm}

In \cite{olli}, Ollivier defines the coarse Ricci curvature between two different points for discrete time Markov chains in the following way:
\begin{defi}[Ollivier]\label{crv} If $P$ is the transition kernel of a Markov chain on a metric space $(E,d)$, the coarse Ricci curvature between $x$ and $y$ is defined by
$$\kappa(x,y)=1-\frac{W_1(\delta_x.P,\delta_y.P)}{d(x,y)}.$$
\end{defi}

A natural generalization of this quantity for continuous-time Markov processes is the following:
\begin{defi}\label{curv} The coarse Ricci curvature between $x$ and $y$ is defined by:
$$\kappa(x,y)=\varliminf_{t\rightarrow0}\frac{1}{t}\left(1-\frac{W_1(P_x^t,P_y^t)}{d(x,y)}\right)$$
where $P_x^t=\delta_x.P^t$.
We also denote
$$\bar\kappa(x,y)=\varlimsup_{t\rightarrow0}\frac{1}{t}\left(1-\frac{W_1(P_x^t,P_y^t)}{d(x,y)}\right).$$
\end{defi}

\begin{rmk}\label{sad}
For every $(x,y,z)\in E^3$, we have the inequality
$$\kappa(x,z)\geq\frac{d(x,y)\kappa(x,y)+d(y,z)\kappa(y,z)}{d(x,z)}$$
(this trivially comes from the triangular inequality for $W_1$).
This property is not always satisfied by $\bar\kappa$ as we will see in the example below.
This is the reason why we choose the liminf in the definition of Ricci curvature.

This inequality is particularly interesting when $d(x,z)=d(x,y)+d(y,z)$ because in this case, the right-hand term of the inequality is a convex combination of $\kappa(x,y)$ and $\kappa(y,z)$, so we have $\kappa(x,z)\geq \min(\kappa(x,y),\kappa(y,z))$.

So in the case of $\varepsilon$-geodesic spaces (see Proposition 19 in \cite{olli}), the infimum of $\kappa(x,y)$ on $E^2$ equals the infimum of $\kappa(x,y)$ for the couples $(x,y)$ such that $d(x,y)\leq\varepsilon$.
So we only have to pay attention to "local" curvature.
\end{rmk}

\begin{xpl}
Let us illustrate the difference between $\kappa$ and $\bar\kappa$.
Let $f:\Rl\mapsto\Rl$ be an increasing continuous function.
Then the deterministic kernel defined on the space $f(\Rl)$ by $P_x^t=\delta_{f(f^{-1}(x)+t)}$ is Markovian.
We choose $f$ such that there exist $t_1<t_2<t_3$ such that
$$f(t)=\left\{\begin{array}{ll}f(t_1)+(t-t_1)(1+\frac{1}{\sqrt{2}}\sin(\ln(|t-t_1|)))&\textrm{in a neighborhood of }t_1\\
f(t_2)+t-t_2&\textrm{in a neighborhood of }t_2\\
f(t_3)+(t-t_3)(1+\frac{1}{\sqrt{2}}\sin(\ln(|t-t_3|)))&\textrm{in a neighborhood of }t_3
\end{array}\right.$$
\begin{center}
\includegraphics{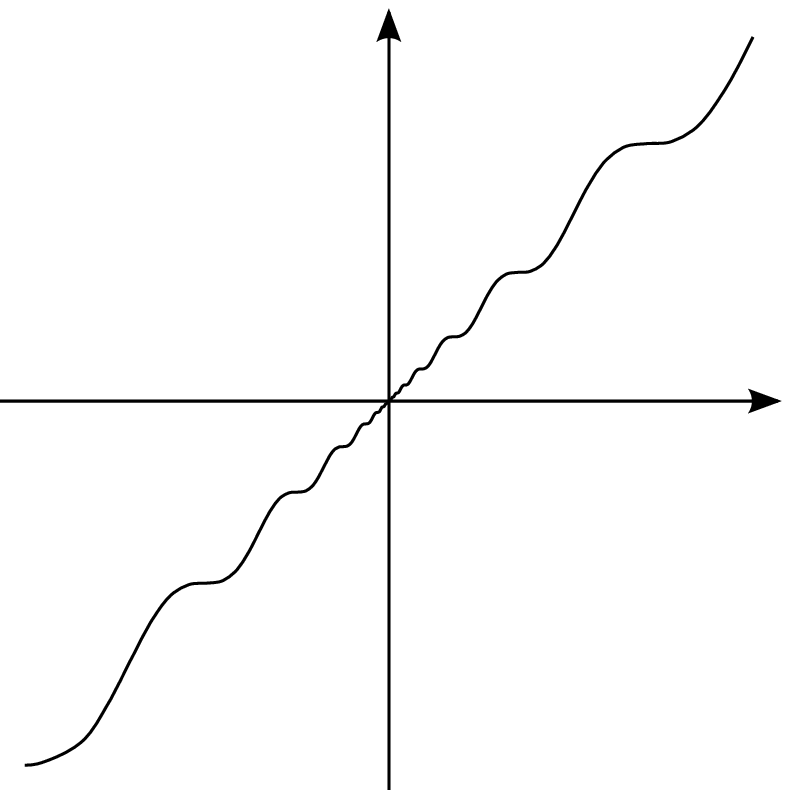}
\end{center}
The graph above is the one of the function $x\mapsto x(1+\frac{\sin(\frac{2\pi\ln(|x|)}{\ln(2)})}{\sqrt{1+\frac{4\pi^2}{\ln^2(2)}}})$, it illustrates the behaviour of $f$ on the neighborhoods of $t_1$ and $t_3$.

Then, if we note $x=f(t_1)$, $y=f(t_2)$, $z=f(t_3)$, the curvatures are $\kappa(x,y)=-\frac{1}{\sqrt{2}(y-x)}$, $\kappa(y,z)=-\frac{1}{\sqrt{2}(z-y)}$ and $\kappa(x,z)=0$, whereas $\bar\kappa(x,y)=\frac{1}{\sqrt{2}(y-x)}$, $\bar\kappa(y,z)=\frac{1}{\sqrt{2}(z-y)}$ and $\bar\kappa(x,z)=0$.
We have $\kappa(x,z)\geq\inf(\kappa(x,y),\kappa(y,z))$, as stated in Remark \ref{sad}, and the same is not true for $\bar\kappa$.
\end{xpl}

\section{$W_1$ contraction in positive coarse Ricci curvature}

It is known (\cite{bd97,cw94,dob70,dob96,olli}) that for Markov Chains, a positive coarse Ricci curvature implies that the Markov operator acting on measures is contractive for the $W^1$ distance.
The following Theorem is a generalization of this result to continuous-time Markov processes:

\begin{thm}\label{inf2} Let $P^t$ be the semigroup of a left-continuous Markov process satisfying $\bar{\kappa}(x,y)\geq K>-\infty$ for all $(x,y)$ in $E^2$ with $x\neq y$.
Then we have: 
$$\forall(x,y)\in E^2,W_1(P_x^t,P_y^t)\leq d(x,y)\nep^{-Kt}.$$
\end{thm}

The hypothesis of this Theorem states that for every $\varepsilon>0$ and $(x,y)\in E^2$, there exist $t<\varepsilon$ such that $W_1(P_x^t,P_y^t)\leq d(x,y)\nep^{-(K-\varepsilon)t}$, but we do not control how this $t$ depends on $x$ and $y$.
The infimum of this $t$ on every neighborhood of every pair of points could be $0$, so we have to refine Ollivier's proof for the discrete time case.

The hypothesis of this Theorem may seem difficult to check on concrete examples, but we can hope to compute $\kappa$, or $\bar{\kappa}$ thanks to the generator of the Markov process in classical cases, under some assumption about the growth in $t$ of the first momentum of $P_x^t$, as we did in \cite{veys} for diffusion processes on manifolds.

\begin{rmk} In particular, the same inequality $W_1(P_x^t,P_y^t)\leq d(x,y)\nep^{-Kt}$ holds if $\kappa(x,y)\geq K$.
\end{rmk}

\begin{cor}
We have $\inf_{x,y}\bar{\kappa}(x,y)=\inf_{x,y}\kappa(x,y)$ for left-continuous Markov processes.
\end{cor}

\noindent\textbf{Proof of the Corollary :} Since $\bar{\kappa}\geq\kappa$, we trivially have $\inf_{x,y}\bar{\kappa}(x,y)\geq\inf_{x,y}\kappa(x,y)$.
Now set $K=\inf_{x,y}\bar{\kappa}(x,y)$.
Theorem \ref{inf2} tells us that for any $(x,y)$, $W_1(P_x^t,P_y^t)\leq d(x,y)\nep^{-Kt}$, so the definition of $\kappa(x,y)$ implies $$\kappa(x,y)\geq\varliminf_{t\rightarrow 0}\frac{1}{t}\left(1-\frac{d(x,y)\nep^{-Kt}}{d(x,y)}\right)=\varliminf_{t\rightarrow 0}\frac{1-\nep^{-Kt}}{t}=K.$$
Thus $\inf_{x,y}\bar{\kappa}(x,y)\leq\inf_{x,y}\kappa(x,y)$.$\square$

\begin{rmk} Usually in the literature, the processes are chosen right-continuous, but Theorem \ref{inf2} also works when the process admits a left-continuous modification. Indeed, the conclusion of the theorem only depends on the law of the process.
\end{rmk}

So Theorem \ref{inf2} does apply to diffusion processes and to minimal jump processes as defined in \cite{chen}, when they do not explode in a finite time.
Indeed, such processes admit left-continuous modifications.
In the case of minimal jump processes, we just have to replace the value of the process at the time of the jump with the value of the process just before the jump to make it left-continuous, and this is a modification, because for every $t$, the probability that the process jumps at time $t$ is $0$.

Theorem \ref{inf2} also applies to some jump processes with an infinite number of jumps in a finite time, provided the locations of the jumps tend to one state in $E$, from which the jump process restarts, as in the following example.

\begin{xpl}\label{ex1} Take the process on $E=\{0\}\cup\{2^{-n},n\in\N\}$ defined in the following way: jump from state $2^{-n}$ to state $2^{-(n+1)}$ after a time of exponential law $\mathcal{E}(2^n+1)$. As $\sum_{n=0}^\infty\frac{1}{2^n+1}<\infty$, the sum of the times of the jumps converges almost surely.
After this infinite number of jumps, the process restarts at $0$ and then jumps to $1$ after a time of law $\mathcal{E}(\frac{1}{2})$.
\end{xpl}

In this example, we have $\kappa(0,2^{-n})=\frac{3}{2}$ and $\kappa(2^{-n},2^{-(n+1)})=\frac{1}{2}$, so thanks to Remark \ref{sad}, we have $\inf(\kappa(x,y))=\frac{1}{2}$, so we can use Theorem \ref{inf2}.

\begin{cor} Let $(P^t)$ be the semigroup of a Markov process on a Polish space admitting a left-continuous modification.
Assume that $\kappa(x,y)>K>0$ for all $(x,y)\in E^2$ with $x\neq y$, and that there exist some $x_0\in E$, $t_0>0$ and $M>0$ such that $W_1(\delta_{x_0},P_{x_0}^t)<M$ for every $0\leq t\leq t_0$.

Then the Markov process admits an unique equilibrium probability measure.
This equilibrium measure has a finite first moment.
\end{cor}

\noindent \textbf{Proof of the Corollary: }
We consider the process starting at $x_0$, restricted to times which are integer multiples of $t_0$.
Using the $W_1$ contraction implied by Theorem \ref{inf2} ,we can easily prove by induction that $W_1\left(P_{x_0}^{nt_0},P_{x_0}^{(n+1)t_0}\right)\leq M\nep^{-Knt_0}$.
So the sequence $(P_{x_0}^{nt_0})$ is a Cauchy sequence for $W_1$, and then it converges to a limit $\pi$ in the $W_1$ distance, and $\pi$ admits a finite first moment.
Now if $t$ is not an integer multiple of $t_0$, we have $W_1\left(P_{x_0}^{\left\lfloor\frac{t}{t_0}\right\rfloor t_0},P_{x_0}^t\right)\leq M\nep^{-K\left\lfloor\frac{t}{t_0}\right\rfloor t_0}$, and the right hand term tends to $0$ when $t$ tends to infinity.
Thus the family $(P_{x_0}^t)$ also tends to $\pi$.

Now we have for every $T>0$,
$$W_1(\pi.P^t,\pi)\leq W_1(\pi.P^t,P_{x_0}^{T+t})+W_1(P_{x_0}^{T+t},\pi)\leq \nep^{-Kt}W_1(\pi,P_{x_0}^T)+W_1(P_{x_0}^{T+t},\pi)$$
and the right hand term tends to $0$ when $T$ tends to the infinity, so $W_1(\pi.P^t,\pi)=0$ and thus $\pi$ is invariant.

Since $W_1(P_x^t,P_{x_0}^t)\leq\nep^{-Kt}d(x,x_0)$, $(P_x^t)$ converges to $\pi$ in $W_1$ distance and thus in weak convergence topology for every $x\in E$.
Then $\mu.P^t$ converges weakly to $\pi$ for every probability measure $\mu$, including any invariant one.
Thus $\pi$ is the unique equilibrium probability measure.$\square$

Theorem \ref{inf2} implies Theorem \ref{inf} as follows.

\noindent\textbf{Proof of Theorem \ref{inf}:}
Let $\pi$ be the unique reversible probability measure.
Theorem \ref{inf2} implies that the operator $P^t$ acting on the space $\mathrm{Lipsch}_0(\pi)$ of Lipschitz functions with mean $0$ (with respect to $\pi$) has a norm smaller than $\nep^{-Kt}$.
Under the hypothesis of Theorem \ref{inf}, the $L^2$ norm is controlled by the Lipschitz norm because $\Var_\pi(f)\leq \Mn_\pi[(f(x)-f(x_0))^2]\leq\|f\|^2_{\mathrm{Lipsch}}\Mn_\pi[d(x,x_0)^2]$ (keep in mind that $\Mn_\pi[d(x,x_0)^2]$ is assumed to be finite).
Now, for any self-adjoint operator $S$ on a Hilbert space $H$, for any $x\in H, x\neq0$, we have $\frac{\|Sx\|}{\|x\|}\leq\sqrt{\frac{\|S^2x\|}{\|x\|}}$, because $\frac{\|S^2x\|^2}{\|x\|^2}-\frac{\|Sx\|^4}{\|x\|^4}=\frac{\|S^2x-\frac{\|Sx\|^2}{\|x\|^2}x\|^2}{\|x\|^2}\geq0$.
So by induction, we get $\frac{\|Sx\|}{\|x\|}\leq\left(\frac{\|S^{2^n}x\|}{\|x\|}\right)^{\frac{1}{2^n}}$.
As $\pi$ is reversible, $P^t$ is self-adjoint on $L^2_0(\pi)$, so we use this result with $S=P^t$ and $x=f\in \mathrm{Lipsch}_0(\pi)\subset L^2_0(\pi)$: we get
$$\frac{\|P^tf\|_{L^2}}{\|f\|_{L^2}}\leq\left(\frac{\sqrt{\Mn_\pi[d(x,x_0)^2]}\|f\|_{\mathrm{Lipsch}}\nep^{-2^nKt}}{\|f\|_{L^2}}\right)^{\frac{1}{2^n}}.$$
The right hand term tends to $\nep^{-Kt}$ when $n$ tends to infinity.
So we have shown that $\|P^tf\|_{L^2}\leq\nep^{-Kt}\|f\|_{L^2}$ for any $f\in\mathrm{Lipsch}_0(\pi)$.
The probability measure $\pi$ is regular (see, for example \cite{part}), so indicator functions can be approximated in $L^2(\pi)$ norm by Lipschitz functions, so Lipschitz functions are dense in $L^2(\pi)$.
Thus, $\mathrm{Lipsch}_0(\pi)$ is dense in $L^2_0(\pi)$.
The operator $P^t$ is $1$-Lipschitz on $L^2(\pi)$, so it is continuous on $L^2_0(\pi)$, and then $\|P^tf\|_{L^2}\leq\nep^{-Kt}\|f\|_{L^2}$ for any $f\in L^2_0(\pi)$.$\square$

\begin{xpl}
Consider the Brownian motion on the circle $\Rl/2\pi\Z$, equipped with the "Euclidean" distance: $$d(\theta_1,\theta_2)=2\left|\sin\left(\frac{\theta_2-\theta_1}{2}\right)\right|.$$
This distance is not geodesic, so we have to compute $\kappa(\theta_1,\theta_2)$ for each $(\theta_1,\theta_2)$ such that $\theta_2\neq\theta_1+2k\pi$, not only for those such that $|\theta_1-\theta_2|<\varepsilon.$
The distance is smooth and bounded, so we can use the formula in \cite{veys} to compute the coarse Ricci curvature.
The Taylor expansion of the distance is:
$$d(\theta_1+\varepsilon v,\theta_2+\varepsilon w)=d(\theta_1,\theta_2)\left(1+\frac{w-v}{2\tan\left(\frac{\theta_2-\theta_1}{2}\right)}-\frac{(w-v)^2}{8}+O(\varepsilon^3)\right).$$
So from \cite{veys}, the coarse Ricci curvature is $\kappa(\theta_1,\theta_2)=0-\frac{-\frac{1}{4}-\frac{1}{4}}{2}+\sqrt{\frac{1}{4}\times\frac{1}{4}}=\frac{1}{2}$.

The process has a positive curvature with this non-geodesic distance, and Theorem \ref{inf} gives the right spectral gap for the generator $\frac{1}{2}\frac{\der^2}{\der\theta^2}$.
If we use the geodesic distance, we get nothing because curvature is 0.
\end{xpl}

To prove Theorem \ref{inf2}, we need a generalization of stopping times, which we call weak stopping times:
\begin{defi}\label{ta}
Let $X_t$ be a random process on a probability space $\Omega$, and $\F_t$ be its natural filtration.
Let $\F_\infty=\mathfrak{S}((\F_t)_{t\in\Rl_+})$ be the $\sigma$-algebra generated by all the $\F_t$. 
A random variable $T$ is a \emph{weak stopping time} for $X_t$ if there exists a $\sigma$-algebra $\mathcal{G}$ independent of $\F_\infty$ such that $T$ is a stopping time for the filtration $\mathcal{G}_t=\mathfrak{S}(\mathcal{G},\F_t)$, ie a positive real-valued random variable $T$ such that $\forall t\geq0,\{\omega\in\Omega|T(\omega)\leq t\}\in \mathcal{G}_t$.
\end{defi}
\begin{lem}\label{taf}
Let $T$ be a random variable having the form $T=\varphi(\omega,v)$, where $\omega$ and $v$ are independent and $\omega$ is a left-continuous process, and $\forall v, \forall t, \mathbf{1}_{\varphi(\omega,v)\leq t}$ only depends on $\omega|_{[0,t]}$.
Then $T$ is a weak stopping time for the process $\omega$, with $\mathcal{G}$ the $\sigma$-algebra generated by $v$.

Conversely, let $X_t$ be a left-continuous random process, and $T$ be any weak stopping time for $X_t$.
Let $\omega$ and $u$ be two independent random variables on another probability space, having the law of $X$ and the uniform law on $[0,1]$ (and assume that $\omega$ is left-continuous).
Then there exists a measurable function $\varphi(\omega,u)$ such that $(\omega_t,\varphi(\omega,u))$ has the law of $(X_t,T)$.

\end{lem}

\noindent\textbf{Proof: }
The first part of the Lemma is trivial once we note that the measurable sets which depend on $\omega|_{[0,t]}$ and $v$ are in $\mathfrak{S}(\mathcal{G},\mathcal{F}_t)$, with $\mathcal{G}$ the $\sigma$-algebra of events only depending on $v$ and $\F_t$ the natural filtration of $\omega$.

So let us prove the other part of the Lemma.
Take $Y_t=\Mn[\mathbf{1}_{T<t}|\F_\infty]$.
It is $\F_t$-measurable because $T$ is a weak stopping time, so there exists a measurable function $f_t$ such that $Y_t=f_t(X|_{[0,t]})$.
If $t_1<t_2$, we have $\mathbf{1}_{T<t_1}\leq \mathbf{1}_{T<t_2}$, then $Y_{t_1}\leq Y_{t_2}$ almost surely by taking the conditional expectation with respect to $\F_\infty$.
So for all $\omega$ outside an exceptional set $N$ of mesure $0$, the function $t\mapsto f_t(\omega|_{[0,t]})$ is non-decreasing on the subset of rational times, and bounded by $1$.
We define the events $A_t=\{\omega|\exists\omega'\notin N,\omega'|_{[0,t]}=\omega|_{[0,t]}\}$.
We take
$$f'_t(\omega|_{[0,t]})=\sup_{t'<t,t\in\Q}(f_{t'}(\omega|_{[0,t']})\mathbf{1}_{A_t}+\mathbf{1}_{A_t^c}).$$
Then for all $\omega$, $t\mapsto f'_t(\omega|_{[0,t]})$ is non-decreasing and left continuous.
Furthermore, we can write $f'_t(X|_{[0,t]})=\Mn[\mathbf{1}_{T<t}|\F_\infty]$ by using the fact that $\mathbf{1}_{T<t}$ is the limit of the increasing sequence $\mathbf{1}_{T<t'_i}$ for any increasing rational sequence $t'_i$ converging to $t$, and the monotone convergence theorem.
Here, $t\mapsto f_t(\omega|_{[0,t]})$ is a kind of repartition function of the conditional law of $T$ knowing $\omega$.
We just have to take $\varphi(\omega,u)=\sup\{t\in\bar\Rl_+|f'_t(\omega|_{[0,t]})<u\}$. $\square$

\begin{defi}If $(S,\mathcal{A})$ is a measurable space, a \emph{kernel} on $S$ will be a measurable application from $S$ to the set of probability measures on $S$.

If $k$ is a kernel on $S$ and $\mu$ is a finite measure on $S$, $\mu.k$ is the finite measure defined by $\mu.k(A)=\int k(x)(A)\der\mu(x)$.
If $k_1$ and $k_2$ are two kernels on $S$, $k_1\ast k_2$ is the kernel defined by $k_1\ast k_2(x)(A)=\int k_2(y)(A)\der k_1(x)(y)$.\end{defi}

\noindent\textbf{Proof of Theorem \ref{inf2}: }
Let $t>0$ and $\varepsilon>0$.
We will show that for every $x$ and $y$,
$$W_1(P_x^t,P_y^t)\leq d(x,y) \nep^{-(K-\varepsilon)t}.$$

We denote by $M^{(x)}$ the Markov process starting at point $x\in E$, and $M^{(x)}_t$ its value at time $t$.
We consider the set $\mathcal{K}$ of kernels $k$ on $E^2\times[0,t]$ satisfying:
\begin{itemize}
\item $\forall(x,y,s),\int d(X,Y)\nep^{(K-\varepsilon)S}\der k((x,y,s))(X,Y,S)\leq d(x,y)\nep^{(K-\varepsilon)s}$.\\
\item $\forall(x,y,s), k((x,y,s))(E^2\times[0,s[)=0$ (i.e. $k$ is a time increasing kernel)
\item there exist weak stopping times $T$ and $T'$ for the Markov process starting at $x$ and $y$, depending measurably on $(x,y,s)$, such that for any random variable $(X,Y,S)$ having the law $k(x,y,s)$, we have $(X,S)\sim(M^{(x)}_T,T+s)$ and $(Y,S)\sim(M^{(y)}_{T'},T'+s)$.
\end{itemize}

Let $(x_0,y_0)\in E^2$, and $\mathcal{I}$ be the set $\{k((x_0,y_0,0)),k\in\mathcal{K}\}$. Our goal is to prove that there exists an element $(X,Y,S)$ of $\mathcal{I}$ satisfying $S=t$ almost-surely, because this would provide us a coupling between $M^{(x_0)}_t$ and $M^{(y_0)}_t$ satisfying

$$\Mn[d(X,Y)]\leq d(x_0,y_0)\nep^{-(K-\varepsilon)t}.$$

We will prove that $\mathcal{I}$ is an inductive set for a well-chosen order relation, and that any maximal element of $\mathcal{I}$ (whose existence is guaranteed by Zorn's Lemma) satisfies $S=t$ almost surely.

To do this, we will prove some nice properties of $\mathcal{K}$:

\begin{lem}\label{stb} The set $\mathcal{K}$ is stable under $\ast$, and any sequence $(k_i)_{i\in\N^*}$ of elements of $\mathcal{K}$ satisfies that the sequence of the products $(k_1\ast k_2\ast \dots \ast k_n)_{n\in\N^*}$ has a limit $k_\infty$ in $\mathcal{K}$, in the sense that for all $(x,y,s)\in E^2\times[0,t]$, the sequence $(k_1\ast k_2\ast\dots\ast k_n)((x,y,s))$ weakly converges to $k_\infty((x,y,s))$.
\end{lem}

\noindent\textbf{Proof of the Lemma:}

Let $(k_i)_{i\in\N}$ be a sequence of elements of $\mathcal{K}$.
Let $(x,y,s)\in E^2\times[0,t]$, and $(X_i,Y_i,S_i)_{i\in\N^*}$ be a Markov chain with non-stationary kernel $k_i$.
Then $(X_n,Y_n,S_n)$ has law $(k_1\ast k_2\ast\dots\ast k_n)(x,y,s)$.
The quantity $\Mn[d(X_i,Y_i)\nep^{(K-\varepsilon)S_i}]$ is non-increasing in $i$.
Indeed, we have \begin{align*}\Mn[d(X_{i+1},Y_{i+1})\nep^{(K-\varepsilon)S_{i+1}}]&=\Mn[\Mn[d(X_{i+1},Y_{i+1})\nep^{(K-\varepsilon)S_{i+1}}|(X_i,Y_i,S_i)]]\\
&\leq\Mn[d(X_i,Y_i)\nep^{(K-\varepsilon)S_i}].\end{align*}
Thus $\Mn[d(X_i,Y_i)\nep^{(K-\varepsilon)S_i}]\leq d(x,y)\nep^{(\kappa-\varepsilon)s}$.

Taking $i=2$ in the previous expression just says that $k_1\ast k_2$ satisfies the first condition in the definition of $\mathcal{K}$.
We will prove below that the sequence $(X_i,Y_i,S_i)$ converges almost surely, and we note $(X_\infty,Y_\infty,S_\infty)$ the limit of this sequence.
Then, thanks to the monotone convergence theorem, we will get $\Mn[d(X_\infty,Y_\infty)\nep^{(K-\varepsilon)S_\infty}]\leq d(x,y)\nep^{(K-\varepsilon)s}$, which is the first condition to check for proving that $k_\infty$ belongs to $\mathcal{K}$.

Now we construct variables $(X'_i,Y'_i,S'_i)$ having the same law as $(X_i,Y_i,S_i)$ over the appropriate probability spaces to prove we have the weak stopping times required by the definition of $\mathcal{K}$.

Let $\Omega^0$ be the space of left-continuous functions from $\Rl_+$ to $E$.
Let us apply Lemma \ref{taf} to the weak stopping times $T$ coming from the third condition of the definition of $\mathcal{K}$ applied to $k_i$: there exist measurable functions $\varphi_i$ from $E^2\times[0,t]\times\Omega^0\times[0,1]$ to $\Rl_+$ such that for every $(x',y',s')\in E^2\times[0,t]$, $\mathbf{1}_{\varphi_i(x',y',s',\omega,u)\leq s''}$ does not depend on the values of $\omega$ for times greater than $s''$, and if we choose $\omega$ and $u$ two independent random variables with laws $\Pb_{x'}$ (the law of the Markov process starting at $x'$) and the uniform law on $[0,1]$, then $(\omega(\varphi_i(x',y',s',\omega,u)),\varphi_i(x',y',s',\omega,u))$ has the law of $(X',S'-s')$ where $(X',Y',S')$ has the law $k_i(x',y',s')$.

Using disintegration of measure on $k_i(x',y',s')$ gives us the existence of a conditional law of $Y'$ knowing $X'=x''$ and $S'=s''$ when $(X',Y',S')$ has the law $k_i(x',y',s')$, and this conditional law depends measurably on $(x',y',s',x'',s'')$. 
Furthermore, as $E$ is a Polish space, any probability measure on $E$ is the law of $f(u)$ with $f:[0,1]\mapsto E$ a measurable function, depending measurably on the probability measure on $E$ and $u$ is a random variable with the uniform law on $[0,1]$.
Then there exist measurable functions $\psi_i:E^2\times[0,t]\times E\times[0,t]\times[0,1]\mapsto E$ such that the law of $\psi_i(x',y',s',x'',s'',u)$ with $u$ a uniform random variable on $[0,1]$ is the conditional law of $Y'$ knowing $X'=x'',S'=s''$, where $X',Y',S'$ has the law $k_i(x',y',s')$.

We take $\Omega=\Omega_0\times[0,1]^\N\times[0,1]^\N$, and we will denote by $(\omega,(u_i)_{i\in\N},(v_j)_{j\in\N})$ the typical element of this set.
We put on the space $\Omega$ the probability measure $\Pb_x\otimes\mathcal{U}([0,1])^{\otimes\N}\otimes\mathcal{U}([0,1])^{\otimes\N}$, which depends on $(x,y,s)$ in a measurable way.
Now we define the following random variables over $\Omega$:
$$\begin{array}{l}X'_0=x\\
Y'_0=y\\
S'_0=s\\
\omega_0=\omega\\
X'_{i+1}=\omega(\varphi_{i+1}(X'_i,Y'_i,S'_i,\omega_i,u_i))\\
S'_{i+1}=S_i+\varphi_{i+1}(X'_i,Y'_i,S'_i,\omega_i,u_i)\\
Y'_{i+1}=\psi_{i+1}(X'_i,Y'_i,S'_i,X'_{i+1},S'_{i+1},v_i)\\
\forall s',\omega_{i+1}(s')=\omega(s'+\varphi_{i+1}(X'_i,Y'_i,S'_i,\omega_i,u_i))\end{array}$$
(in other words, $\omega_{i+1}$ is $\omega_i$ "shifted" by $S'_{i+1}-S'_i$).

We prove by induction that for all $n$, $S'_n-s$ is a weak stopping time with $\mathcal{G}=\mathfrak{S}((u_i)_{i\in\N},(v_j)_{j\in\N})$ and the conditional law of $(\omega_n,(u_{n+i}),(v_{n+j}))$ knowing $(X'_0,Y'_0,S'_0,\dots,X'_n,Y'_n,S'_n)$ is $\Pb_{X'_n}\otimes\mathcal{U}([0,1])^{\otimes\N}\otimes\mathcal{U}([0,1])^{\otimes\N}$.

The case $n=0$ is trivial.
If we fix $u_0,\dots,u_n,v_0,\dots,v_n$ and take $s\leq s'\leq t$, we have to show that the function $\mathbf{1}_{S'_{n+1}\leq s'}$ does not depend on the values of $\omega$ for times greater than $s'-s$.
Because of the property of $\varphi_{n+1}$, we know that if $S'_n<s'$, the function $\mathbf{1}_{S'_{n+1}\leq s'}$ does not depend on the values of $\omega_{n}$ for times greater than $s'-S'_n$, that is, on the values of $\omega$ for times greater than $s'-s$.
The induction hypothesis tells that $S'_n-s$ is a weak stopping time with $\mathcal{G}=\mathfrak{S}((u_i)_{i\in\N},(v_j)_{j\in\N})$, and then the event $S'_n>s'$ does not depend on values of $\omega$ for times greater than $s'$.
So $S'_{n+1}$ is a weak stopping time with $\mathcal{G}=\mathfrak{S}((u_i)_{i\in\N},(v_j)_{j\in\N})$.
We can use the Markov property, so the conditional law of $\omega_{n+1}$ knowing $(X'_0,Y'_0,S'_0,\dots,X'_{n+1},Y'_{n+1},S'_{n+1})$ is $\Pb_{X'_{n+1}}$.
As $(X'_k,Y'_k,S'_k)$ only depends on $\omega$ and the $u_i$'s and $v_j$'s with $i$ and $j$ smaller than $k-1$, so the subsequence $(u_{n+1+i},v_{n+1+j})$ is independent of $(X'_0,Y'_0,S'_0,\dots,X'_{n+1},Y'_{n+1},S'_{n+1})$.

So the conditional law of $(X'_{i+1},Y'_{i+1},S'_{i+1})$ knowing $(X'_0,Y'_0,S'_0,\dots,X'_i,Y'_i,S'_i)$ is $k_{i+1}((X'_i,Y'_i,S'_i))$.
Thus $((X'_i,Y'_i,S'_i))_{i\in\N}$ and $((X_i,Y_i,S_i))_{i\in\N}$ have the same law.

The sequence $(S'_i)$ is non-decreasing and bounded by $t$, so it converges almost surely to a limit $S'_\infty$, which is also a weak stopping time, because the supremum of a family of stopping times for the filtration $\mathcal{G}_t$ is a stopping time for the filtration $\mathcal{G}_t$.
Because of the left continuity of $\omega$, the sequence $(X'_i)$ converges to $X'_\infty=\omega(S'_\infty-s)$.

Of course, we can do the same thing by swapping the roles of $x$ and $y$ to define $(X''_i,Y''_i,S''_i)$, and then we have the convergence of $(S''_i)$ to a weak stopping time $S''_\infty$ and the convergence of $(Y''_i)$ to $Y''_\infty=\omega'(S''_\infty-s)$.
So we have proved that $((X_i,Y_i,S_i))$ converges almost surely, thus we have the existence of a limit $k_\infty$ of $k_1\ast k_2\ast\dots\ast k_n$.

The fact that $S'_2$, $S''_2$, $S'_\infty$ and $S''_\infty$ are weak stopping times show us that $k_1\ast k_2$ and $k_\infty$ satisfy the last two points of the definition of $\mathcal{K}$, so they belong to $\mathcal{K}$.$\square$

\noindent\textbf{End of the proof of Theorem \ref{inf2}:}

Let us put the following partial order relation on $\mathcal{I}$: $\mu_1\preceq\mu_2$ if and only if there exists $k\in\mathcal{K}$ so that $\mu_2=\mu_1.k$.
First we check that $\preceq$ is an order relation.
Transitivity of $\preceq$ is due to the fact that $\mathcal{K}$ is stable under $\ast$.
Reflexivity is a consequence that $1_{\ast}\in\mathcal{K}$, with $1_{\ast}:(x,y,s)\rightarrow\delta_{(x,y,s)}$ the trivial kernel.
Antisymmetry is a bit harder to check.
Suppose $\mu\preceq\nu\preceq\mu$.
We have $\nu=k_1.\mu$ and $\mu=k_2.\nu$ with $(k_1,k_2)\in\mathcal{K}^2$.
We construct $(X_0,Y_0,S_0)$ of law $\mu$, $(X_1,Y_1,S_1)$ of conditional law $k_1(X_0,Y_0,S_0)$ knowing $(X_0,Y_0,S_0)$, and $(X_2,Y_2,S_2)$ of conditional law $k_2(X_1,Y_1,S_1)$ knowing $(X_0,Y_0,S_0,X_1,Y_1,S_1)$.
Then we have $S_0\leq S_1\leq S_2$ almost surely, so $\Mn[S_0]\leq\Mn[S_1]\leq\Mn[S_2]=\Mn[S_0]$ ($S_0$ and $S_2$ have the same law).
As $S_1-S_0\geq 0$ almost surely and $\Mn[S_1-S_0]=0$, we have $S_0=S_1$ almost surely.
Since $k_1\in\mathcal{K}$, we then have $k_1(x,y,s)=\delta_{(x,y,s)}$, $\mu$-almost surely, and so $\nu=\mu$.

Now we will prove that $\mathcal{I}$ is an inductive set.
Let $A\subset\mathcal{K}$ be a totally ordered subset.
If $A$ is empty, then $\delta_{(x_0,y_0,0)}\in\mathcal{I}$ is an upper bound of $A$.
Otherwise, we consider $M=\sup_{\mu\in A}\Mn_\mu[S]\in[0,t]$.
If there exists $\mu\in A$ such that $\Mn_\mu[S]=M$, then $\mu$ is the maximum of $A$.
In the remaining case, there exists an increasing sequence $(\mu_i)_{i\in\N}\in A^\N$ such that $\Mn_{\mu_i}[S]\nearrow M$, and so for every $\mu\in A$, there exists $i\in\N$ so that $\mu\preceq\mu_i$, because $A$ is totally ordered and $\mu\preceq\mu'\Rightarrow\Mn_\mu[S]\leq\Mn_{\mu'}[S]$.
Any upper bound of all the $\mu_i$ will be an upper bound of $A$.
For each $i$, there exists $k_i\in\mathcal{K}$ so that $\mu_{i+1}=\mu_i.k_i$.
Then by lemma \ref{stb} $\mu_\infty=\lim_{i\rightarrow\infty}\mu_i$ exists, belongs to $\mathcal{I}$ and we have for each $i$, $\mu_\infty=\mu_i.(k_i\ast k_{i+1}\ast\dots)$.
So $\mu_\infty$ is an upper bound of $A$.

We can apply Zorn's lemma to $\mathcal{I}$ to get a maximal element $\mu_\mathrm{max}$.
Then we have $\mu_\mathrm{max}.k=\mu_\mathrm{max}$ for every $k\in\mathcal{K}$.
Let us prove that under $\mu_{\mathrm{max}}$, $s=t$ almost surely.
To do so, we will construct a particular $k\in\mathcal{K}$ such that for all $(x,y,s)$, we have $s=t$ or $\Pb_{(x',y',s')\sim k(x,y,s)}(s'>s)=1$, and then the fact that $\mu_\mathrm{max}.k=\mu_\mathrm{max}$ implies that $s=t$ almost surely under $\mu_\mathrm{max}$.

By definition of $\bar{\kappa}$, for each $(x,y,s)$ with $s<t$, there exists $0<\eta(x,y,s)\leq t-s$ such that $W_1(\Pb_x^{\eta(x,y,s)},\Pb_y^{\eta(x,y,s)})\leq d(x,y)\nep^{-(K-\varepsilon)\eta(x,y,s)}$ (because $\bar\kappa(x,y)>K-\varepsilon$).
So we have a coupling $\xi(x,y,s)$ between $\Pb_x^{\eta(x,y,s)}$ and $\Pb_y^{\eta(x,y,s)}$ such that $\Mn_{\xi(x,y,s)}[d(X,Y)]\leq d(x,y)\nep^{-(K-\varepsilon)\eta(x,y,s)}$.

It remains to prove that we can choose $\eta(x,y,s)$ and $\xi(x,y,s)$ in a measurable way to get our $k$.
A simple choice for $\eta(x,y,s)$ is the maximal one
$$\sup(\{\eta\in]0,t-s]|W_1(\Pb_x^\eta,\Pb_y^\eta)\leq d(x,y)\nep^{-(K-\varepsilon)\eta}\}),$$
which is measurable.
The fact that this supremum is actually a maximum is due to the existence of a left continuous modification.
Indeed, let $(\eta_i)_{i\in\N}$ be a maximizing sequence for the expression above, and $\eta$ be the supremum.
Let $f$ be any bounded $1$-lipschitz function from $E$ to $\Rl$.
Because of the left continuous modification, $M^{(x)}_{\eta_i}$ converges to $M^{(x)}_\eta$ and $M^{(y)}_{\eta_i}$ converges to $M^{(y)}_\eta$.
So by the dominated convergence theorem, $\Mn[f(M^{(x)}_{\eta_i})]$ converges to $\Mn[f(M^{(x)}_\eta)]$ and $\Mn[f(M^{(y)}_{\eta_i})]$ converges to $\Mn[f(M^{(y)}_\eta)]$.
Thus $\int f\der(\Pb_x^{\eta_i}-\Pb_y^{\eta_i})$ converges to $\int f\der(\Pb_x^\eta-\Pb_y^\eta)$ and this latter is smaller than $d(x,y)\nep^{-(K-\varepsilon)\eta}$.
Then, there exists a measurable way to choose an optimal coupling between two probability measures (Corollary 5.22 in \cite{vill}), thus we can get a measurable $\xi(x,y,s)$.

We can then set $k(x,y,s)=\xi(x,y,s)\otimes\delta_{s+\eta(x,y,s)}$ for $s<t$, and $k(x,y,t)=\delta_{(x,y,t)}$ (because $\eta(x,y,s)\geq0$ is trivially a weak stopping time).
Since $\eta(x,y,s)>0$, and $\mu_\mathrm{max}.k=\mu_\mathrm{max}$, we have $s=t$ $\mu_\mathrm{max}$-almost surely, so $\mu_\mathrm{max}$ provides a coupling between $\Pb_{x_0}^t$ and $\Pb_{y_0}^t$ which satisfies $\Mn[d(X,Y)]\leq d(x_0,y_0)\nep^{-(\kappa-\varepsilon)t}$, so $W_1(\Pb_{x_0}^t,\Pb_{y_0}^t)\leq d(x_0,y_0)\nep^{-(\kappa-\varepsilon)t}$ as needed.
Letting $\varepsilon$ decrease to $0$ gives the conclusion of Theorem \ref{inf2}.$\square$

\end{document}